\documentclass{article}
\usepackage{hyperref}
\usepackage{natbib}
\usepackage{graphicx} 
\usepackage{amsmath, amsthm, amssymb}
\usepackage{xcolor}
\usepackage{multirow}

\usepackage[a4paper, top=2.5cm, bottom=3cm, left=2.5cm, right=2.5cm]{geometry}
\newtheorem{theorem}{Theorem}

\title{New Sphere Packings from the Antipode Construction}
\author{Ruitao Chen\footnote{Alphabetical order. $^\dagger$Correspondence to \texttt{wanglw@pku.edu.cn}.}\quad Jiachen Hu$^*$\quad  Binghui Li$^*$\quad  Liwei Wang$^{*,\dagger}$\quad  Tianyi Wu$^*$
\\
\quad
\\
Peking University}
\date{}

\begin{document}

\maketitle

\begin{abstract}
    In this note, we construct non-lattice sphere packings in dimensions $19$, $20$, $21$, $23$, $44$, $45$, and $47$, demonstrating record densities that surpass all previously documented results in these dimensions. The construction involves applying the antipode method to suboptimal cross-sections of $\Lambda_{24}$ and $P_{48p}$ respectively in those dimensions. 
\end{abstract}

\section{Introduction}

The sphere packing problem in $n$-dimensional Euclidean space \(\mathbb{R}^n\) asks for an arrangement of non-overlapping congruent spheres that fills the largest possible fraction of space. For a packing \(\mathcal{P}\) of spheres of radius \(r > 0\) centered at points in a set \(\mathcal{C} \subset \mathbb{R}^n\), its density is defined as  
\[
\Delta(\mathcal{P}) = \limsup_{R \to \infty} \frac{\mathrm{Vol}(B_R \cap \bigcup_{x \in \mathcal{C}} B_r(x))}{\mathrm{Vol}(B_R)},
\]  
where \(B_R\) is the ball of radius \(R\) centered at the origin, \(B_r(x)\) is the ball of radius \(r\) centered at the point $x$, and $\mathrm{Vol}(\cdot)$ denotes the volume. Additionally, the center density is given by  
\[
\delta(\mathcal{P}) = \Delta(\mathcal{P})/V_n,
\]  
where $V_n$ is the volume of the unit ball in $\mathbb{R}^n$. The center density provides a useful measure of the distribution of sphere centers in space. The optimal center density in dimension \(n\), denoted \(\delta_n\), is the supremum of \(\delta(\mathcal{P})\) over all such packings in \(\mathbb{R}^n\).  

Exact values of \(\delta_n\) are known only in a few dimensions. In \(\mathbb{R}^1\), the optimal packing simply places intervals end-to-end, yielding center density \(\delta_1 = 2^{-1}\). In \(\mathbb{R}^2\), the hexagonal lattice achieves the optimal density \(\delta_2 = 2^{-1}\cdot 3^{-1/2}\), and in \(\mathbb{R}^3\), the face-centered cubic packing attains the maximal density \(\delta_3 = 2^{-5/2}\), as proved by Hales \citep{hales2005proof} for the Kepler conjecture. 
It is shown that in \(\mathbb{R}^8\), the \(\Lambda_8\) lattice achieves the optimal density $\delta_{8} = 2^{-4}$ \citep{viazovska2017sphere}, and in \(\mathbb{R}^{24}\), the Leech lattice $\Lambda_{24}$ is optimal with $\delta_{24}=1$ \citep{cohn2017sphere}. The optimality of these results was established using techniques from Fourier analysis and modular forms, marking breakthroughs in discrete geometry. 

While exact values of \(\delta_n\) remain unknown in most dimensions, substantial progress has been made in constructing dense packings in low dimensions. The densest lattice packings known in most dimensions below $24$ are laminated lattices $\Lambda_n$ or lattices $K_n$ \citep{leech1964some,leech1967notes}. All these lattices are cross-sections of $\Lambda_{24}$. Some dense non-lattice packings are constructed based on nonlinear binary codes, setting records in dimension $10$, $11$, and $13$ \citep{leech1970new,leech1971sphere}. Construction $\mathrm{B}^*$, a variant of Construction $\mathrm{B}$ of \citet{conway1999sphere} based on two orthogonal binary codes, yields non-lattice packings in dimension $18$ \citep{bierbrauer2000dense}, $20$ \citep{vardy1995new}, and $27$–$30$ \citep{vardy1999density}. The antipode construction, introduced by \citet{conway1996antipode}, can be used to generate a packing composed of translated copies of a base lattice. In dimensions $20$, $22$, and $44$–$47$, the non-lattice packings generated via the antipode construction are denser than their respective base lattices. 

In this note, we present denser non-lattice sphere packings in dimensions $19$, $20$, $21$, $23$, $44$, $45$, and $47$, which surpass the best previously known. Our construction applies the antipode method to carefully selected suboptimal cross-sections of self-dual lattices \(\Lambda_{24}\) \citep{leech1967notes} and \(P_{48p}\) \citep{conway1999sphere}. 
Our results are briefly summarized in Table \ref{tab: results}.

\begin{table}[tbh]
    \centering
    \begin{tabular}{|c|c|c|c|c|}
    \hline
    \multirow{2}{*}{Dimension} & \multicolumn{2}{c|}{Previously best known results} & \multicolumn{2}{c|}{Our results}   \\
    \cline{2-5}
    & Construction & Center density & Construction & Center density
     \\
     \hline 
     19 & Laminated lattice & 0.08839 & Antipode & 0.08896 \\
     \hline 
     20 & Construction $\mathrm{B}^*$/Antipode & 0.13154 & Antipode & 0.15593 \\
     \hline 
     21 & Laminated lattice & 0.17678 & Antipode & 0.21004 \\
     \hline 
     23 & Laminated lattice & 0.50000 & Antipode & 0.50049 \\
     \hline 
     44 & Antipode & 472.799 & Antipode & 509.619 \\
     \hline 
     45 & Antipode & 974.700 & Antipode & 1243.46 \\
     \hline 
     47 & Antipode & 5788.81 & Antipode & 5925.98 \\
     \hline
    \end{tabular}
    \caption{A brief summary of our results.}
    \label{tab: results}
\end{table}

\section{The Antipode Construction}

The antipode construction, introduced by \citet{conway1996antipode}, is a method for generating dense sphere packings through lattice translates. The underlying principle of this construction involves considering the union of multiple translated instances of a base lattice, which serve as the centers for the spheres. Typically, the base lattice is selected as the cross-section of a dense, self-dual lattice (such as the Leech lattice $\Lambda_{24}$ and the lattice $P_{48p}$). 

Let \( \Lambda \) be an $n$-dimensional lattice with minimal norm \( \mu \). The antipode construction begins with an orthogonal decomposition of the Euclidean space \( \mathbb{R}^n \) into two subspaces, \( k \)-dimensional subspace \( U \) and \( l \)-dimensional subspace \( V \) (\( n = k + l \) and \( \mathbb{R}^n = U \oplus V \)), such that \( K = \Lambda \cap U \) and \( L = \Lambda \cap V \) are $k-$dimensional and $\ell-$dimensional cross-sections of the lattice \( \Lambda \) respectively. We define the lattices $M:=\pi_U(\Lambda)$ and $N:=\pi_V(\Lambda)$, where $\pi_U$ and $\pi_V$ denote the projection operators from $\mathbb{R}^{n}$ onto the subspaces $U$ and $V$, respectively. The construction proceeds by finding a subset \( S = \{u_1, u_2, \dots, u_s\} \) (where $s:= |S|$ is the size of the set $S$) of the lattice \( M \), such that the squared Euclidean distance between any two distinct points in \( S \) is bounded above by a constant \( \beta < \mu\). Finally, the $l$-dimensional sphere packing $\mathcal{A}(S)$ in the space $V$ is obtained by constructing spheres of radius $\sqrt{\mu - \beta}/2$ centered at points in the set $\left\{\pi_V(w): w \in \Lambda, \pi_U(w) \in S\right\}$, which is the union of $s$ translated copies of the lattice $L$. 

\begin{theorem}
\label{th: main}
    The center density of the $l$-dimensional antipode packing $\mathcal{A}(S)$ is given by
    \begin{equation}
        \delta(\mathcal{A}(S)) = s\sqrt{\frac{\det M}{\det \Lambda}}\left(\frac{\mu-\beta}{4}\right)^{l/2}. \nonumber
    \end{equation}
\end{theorem}

The method allows the packing to achieve higher densities than the original lattice packing by leveraging points in \( S \). The success of the construction depends on choosing an appropriate set \( S \) that balances the number of translates $s$ and the maximum distance $\beta$ between points in $S$. Typically, we set $s=k+1$ in this note. 


\section{Results}

Previously best-known results in dimensions $20$, $22$, and $44$–$47$ obtained by the antipode method are translations of the densest cross-sections of $\Lambda_{24}$ and $P_{48p}$. However, we found that translations of suboptimal cross-sections yield denser packing in some dimensions.

For dimensions $23$, $21$, $20$, and $19$, we focus on $\Lambda=\Lambda_{24}$, the Leech lattice in $\mathbb{R}^{24}$. The generator matrix of $\Lambda_{24}$ in standard MOG coordinates (Figure $4.12$ of \citet{conway1999sphere}) is used. For dimensions $47$, $45$, and $44$, we use $\Lambda=P_{48p}$, a self-dual lattice with minimal norm $6$. When representing vectors, we use the abbreviation $m^n$ to denote $n$ consecutive $m$. 

In what follows, we present a special case of the general antipode construction (Theorem \ref{th: main}). Throughout this note, all subsequent constructions are derived from Theorem \ref{th: simple}.
\begin{theorem}
Let $\Lambda$ denote an $n$-dimensional self-dual lattice with minimal norm $\mu$. Assume that the lattice $K$ is a $k$-dimensional cross-section of $\Lambda$, characterized by a generator matrix $\mathbf{A}$ and Gram matrix $\mathbf{K}$. Define $\mathbf{M} = \mathbf{K}^{-1}$. If there exists a constant $0<\beta<\mu$ satisfying
\begin{align*}
\mathbf{M}_{i,i}\leq \beta,1\leq i\leq k, \quad\mathrm{and}\quad \mathbf{M}_{i,i}+\mathbf{M}_{j,j}-2\mathbf{M}_{i,j}\leq \beta,1\leq i<j\leq k.
\end{align*}
Then we can construct an $(n-k)$-dimensional sphere packing of center density
\begin{align*}
\delta = (k+1)\sqrt{\det\mathbf{M}}\left(\frac{\mu-\beta}{4}\right)^{(n-k)/2}.
\end{align*}
\label{th: simple}
\end{theorem}

\begin{proof}
Since $\Lambda$ is self-dual, we have lattice $M=K^*$. Note that $\mathbf{M}$ is the Gram matrix of the lattice $M$. We assume $v_1, v_2, \dots,v_k$ are rows of the generator matrix corresponding to $\mathbf{M}$. Let $v_0$ be the origin and $S=\{v_0,v_1,v_2,\dots,v_k\}$. Then 
\begin{align*}
d^2(v_0,v_i)=\mathbf{M}_{i,i}\leq \beta,1\leq i\leq k, \quad\mathrm{and}\quad d^2(v_i,v_j)=\mathbf{M}_{i,i}+\mathbf{M}_{j,j}-2\mathbf{M}_{i,j}\leq \beta,1\leq i<j\leq k.
\end{align*}
The theorem is then proved by Theorem \ref{th: main}.
\end{proof}

For dimensions $23$, $21$, $20$, and $19$, we take cross-sections of $\Lambda_{24}$ as the base lattices. We provide the generator matrix $\mathbf{A}$ and the matrix $\mathbf{M}$ in each dimension. 

\paragraph{Dimension $23$.} We take $\mathbf{K}=[6]$ and $\beta = 1/6$, then get a $23$-dimensional packing with center density 
\begin{align*}
\delta = 23^{11.5}\cdot2^{-34}\cdot3^{-12} = 0.50049\ldots\,.
\end{align*}

\paragraph{Dimension $21$.} We let $\beta=5/12$ and
\begin{align*}
\mathbf{A} = \frac{1}{\sqrt{2}}
\begin{bmatrix}
2&0&0^6&2&0^{15}\\
0&2&0^6&-2&0^{15}\\
-1&-1&1^{6}&0&0^{15}
\end{bmatrix}
,\quad \mathbf{M}=\frac{1}{12}
\begin{bmatrix}
5& 3& 2\\
3& 5& 2\\
2& 2& 4\\
\end{bmatrix}.
\end{align*}
Since $\det\mathbf{M}=1/36$, we obtain a $21$-dimensional sphere packing with center density 
\begin{align*}
\delta = 43^{10.5}\cdot2^{-41}\cdot3^{-11.5} = 0.21004\ldots\,.
\end{align*}

\paragraph{Dimension $20$.} For simplicity, we use $+$ and $-$ to denote $1$ and $-1$, respectively, in the matrix $\mathbf{A}$. We take $\beta=2/5$ and 
\begin{align*}
\mathbf{A} = \frac{1}{\sqrt{2}}
\begingroup
\setlength{\arraycolsep}{3pt}
\begin{bmatrix}
++&++&-+&-+&00&00&00&00&0^8\\
-+&--&00&00&++&+-&00&00&0^8\\
+-&00&+-&00&+-&00&+-&00&0^8\\
-0&+0&+0&+0&-0&-0&-0&+0&0^8\\
\end{bmatrix}
\endgroup
, \quad\mathbf{M}=\frac{1}{5}
\begin{bmatrix}
2& 1& 1& 1\\
1& 2& 1& 1\\
1& 1& 2& 1\\
1& 1& 1& 2
\end{bmatrix}.
\end{align*}
Since $\det\mathbf{M}=1/125$, we obtain a $20$-dimensional sphere packing with center density
\begin{align*}
\delta = 3^{20}\cdot2^{-10}\cdot5^{-10.5} = 0.15593\ldots\,.
\end{align*}

\paragraph{Dimension $19$.} We take $\beta=8/15$ and 
\begin{align*}
\mathbf{A} = \frac{1}{\sqrt{2}}
\begingroup
\setlength{\arraycolsep}{3pt}
\begin{bmatrix}
0&0&0^2&2&0&0^2&2&0&0^2&0^{12}\\
2&0&0^2&0&2&0^2&0&0&0^2&0^{12}\\
0&2&0^2&0&0&0^2&0&2&0^2&0^{12}\\
-1&-1&-1^2&-1&-1&1^2&0&0&0^2&0^{12}\\
-1&-1&1^2&0&0&0^2&-1&-1&1^2&0^{12}
\end{bmatrix}
\endgroup
,\quad \mathbf{M}=\frac{1}{60}
\begin{bmatrix}
21& 9& 9& 12& 12\\
9& 29& 13& 20& 16\\
9& 13& 29& 16& 20\\
12& 20& 16& 32& 16\\
12& 16& 20& 16& 32\\
\end{bmatrix}.
\end{align*}
Since $\det\mathbf{M}=1/300$, we can construct a $19$-dimensional sphere packing with center density
\begin{align*}
\delta = 13^{9.5}\cdot3^{-9}\cdot5^{-10.5} = 0.08896\ldots\,.
\end{align*}

For dimensions $47$, $45$, and $44$, we take cross-sections of $P_{48p}$ as base lattices. For simplicity, we only state the Gram matrix $\mathbf{K}$ and its corresponding inverse $\mathbf{M}$.

\paragraph{Dimension $47$.} We take $\mathbf{K}=[8]$ and $\beta = 1/8$, then construct a $47$-dimensional packing with center density 
\begin{align*}
\delta = 47^{23.5}\cdot2^{-118}=5925.98\ldots\,.
\end{align*}

\paragraph{Dimension $45$.} We take $\beta=1/4$ and
\begin{align*}
\mathbf{K} = 
\begin{bmatrix}
6&-2&-2\\
-2&6&-2\\
-2&-2&6
\end{bmatrix}
,\quad \mathbf{M}=\frac{1}{8}
\begin{bmatrix}
2&1&1\\
1&2&1\\
1&1&2
\end{bmatrix}.
\end{align*}
Since $\det\mathbf{M}=1/128$, we obtain a $45$-dimensional sphere packing with center density
\begin{align*}
\delta = 23^{22.5}\cdot2^{-91.5} = 1243.46\ldots\,.
\end{align*}

\paragraph{Dimension $44$.} We take $\beta=16/55$ and 
\begin{align*}
\mathbf{K} =
\begin{bmatrix}
6&-2&-1&-1\\
-2&6&-2&-1\\
-1&-2&6&-2\\
-1&-1&-2&6
\end{bmatrix}
,\quad \mathbf{M}=\frac{1}{55}
\begin{bmatrix}
14&8&7&6\\
8&16&9&7\\
7&9&16&8\\
6&7&8&14
\end{bmatrix}.
\end{align*}
Since $\det\mathbf{M}=1/605$, we can construct a $44$-dimensional sphere packing with center density
\begin{align*}
\delta = 157^{22}\cdot2^{-22}\cdot5^{-21.5}\cdot11^{-23} = 509.619\ldots\,.
\end{align*}


\paragraph{Remark.} For dimension $18$, we found a non-lattice packing with density $5^{9}\cdot2^{-8}\cdot3^{-10.5}=0.07460\ldots$ using antipode construction, which is denser than $\Lambda_{18}$ but less dense than $\mathcal B^*_{18}$. Denser sphere packings in dimensions below $44$ are likely to be constructed using the proposed method.

\bibliographystyle{ims}
\bibliography{ref}

\begin{thebibliography}{12}
\expandafter\ifx\csname natexlab\endcsname\relax\def\natexlab#1{#1}\fi
\expandafter\ifx\csname url\endcsname\relax
  \def\url#1{\texttt{#1}}\fi
\expandafter\ifx\csname urlprefix\endcsname\relax\def\urlprefix{}\fi

\bibitem[{Bierbrauer and Edel(2000)}]{bierbrauer2000dense}
\text{Bierbrauer, J.} and \text{Edel, Y.} (2000).
\newblock Dense sphere packings from new codes.
\newblock \textit{Journal of Algebraic Combinatorics}, \textbf{11} 95--100.

\bibitem[{Cohn et~al.(2017)Cohn, Kumar, Miller, Radchenko and Viazovska}]{cohn2017sphere}
\text{Cohn, H.}, \text{Kumar, A.}, \text{Miller, S.}, \text{Radchenko, D.} and \text{Viazovska, M.} (2017).
\newblock The sphere packing problem in dimension 24.
\newblock \textit{Annals of mathematics}, \textbf{185} 1017--1033.

\bibitem[{Conway and Sloane(1996)}]{conway1996antipode}
\text{Conway, J.} and \text{Sloane, N.} (1996).
\newblock The antipode construction for sphere packings.
\newblock \textit{Inventiones mathematicae}, \textbf{123} 309--313.

\bibitem[{Conway and Sloane(1999)}]{conway1999sphere}
\text{Conway, J.~H.} and \text{Sloane, N. J.~A.} (1999).
\newblock \textit{Sphere Packings, Lattices and Groups}, vol. 290 of \textit{Grundlehren der mathematischen Wissenschaften}.
\newblock 3rd ed. Springer-Verlag, New York.

\bibitem[{Hales(2005)}]{hales2005proof}
\text{Hales, T.~C.} (2005).
\newblock A proof of the kepler conjecture.
\newblock \textit{Annals of mathematics} 1065--1185.

\bibitem[{Leech(1964)}]{leech1964some}
\text{Leech, J.} (1964).
\newblock Some sphere packings in higher space.
\newblock \textit{Canadian Journal of Mathematics}, \textbf{16} 657--682.

\bibitem[{Leech(1967)}]{leech1967notes}
\text{Leech, J.} (1967).
\newblock Notes on sphere packings.
\newblock \textit{Canadian Journal of Mathematics}, \textbf{19} 251--267.

\bibitem[{Leech and Sloane(1970)}]{leech1970new}
\text{Leech, J.} and \text{Sloane, N.} (1970).
\newblock New sphere packings in dimensions 9--15.
\newblock \textit{Bulletin of the American Mathematical Society}, \textbf{76} 1006--1010.

\bibitem[{Leech and Sloane(1971)}]{leech1971sphere}
\text{Leech, J.} and \text{Sloane, N.} (1971).
\newblock Sphere packings and error-correcting codes.
\newblock \textit{Canadian Journal of Mathematics}, \textbf{23} 718--745.

\bibitem[{Vardy(1995)}]{vardy1995new}
\text{Vardy, A.} (1995).
\newblock A new sphere packing in 20 dimensions.
\newblock \textit{Inventiones mathematicae}, \textbf{121} 119--133.

\bibitem[{Vardy(1999)}]{vardy1999density}
\text{Vardy, A.} (1999).
\newblock Density doubling, double-circulants, and new sphere packings.
\newblock \textit{Transactions of the American Mathematical Society}, \textbf{351} 271--283.

\bibitem[{Viazovska(2017)}]{viazovska2017sphere}
\text{Viazovska, M.~S.} (2017).
\newblock The sphere packing problem in dimension 8.
\newblock \textit{Annals of mathematics} 991--1015.

\end{thebibliography}

\end{document}